\documentclass[11pt]{amsart}
\usepackage{mathrsfs}
\usepackage{epic,eepic,epsf,epsfig}
\usepackage{amsfonts,srcltx,mathrsfs}

\usepackage{amsmath}
\usepackage{amssymb}
\usepackage{amsbsy}
\usepackage{float}
\usepackage{graphicx}
\usepackage{amsfonts}
\usepackage{color}

\usepackage{tikz}

\usepackage{url}
\newtheorem{lem}{Lemma}[section]
\newtheorem{theorem}[lem]{Theorem}

\newtheorem{cor}[lem]{Corollary}

\newtheorem{prob}[lem]{Problem}

\newtheorem{case}{Case}
\newtheorem{prop}[lem]{Proposition}

\def\a{\alpha}    
\def\r{\rho} \def\s{\sigma}

\def\olg{\overline G}

\def\di{\bigm|} \def\lg{\langle} \def\rg{\rangle}
\def\nd{\mathrel{\bigm|\kern-.7em/}}

\def\f{\noindent}

\def\PSL{\hbox{\rm PSL}}

\def\PGL{\hbox{\rm PGL}}
\def\PGaL{\hbox{\rm P\mbox{$\Gamma$}L}}
\def\G{\hbox{\rm G}}
\def\Sz{\hbox{\rm Sz}}

\def\Aut{\hbox{\rm Aut}}

\def\Aut{\hbox{\rm Aut}}

\def\mod{\hbox{\rm mod }}
\def\Core{\hbox{\rm Core}}

\def\PGL{\hbox{\rm PGL}}

\def\demo{\f {\bf Proof.}\hskip10pt}

\def\mz{{\mathbb Z}}
\def\P{\mathcal {P}}

\begin{document}

\title{Regular $3$-polytopes of order $2^np$}

\author{Dong-Dong Hou}
\address{Dong-Dong Hou, Department of Mathematics, Shanxi Normal University, TaiYuan,
100044, P.R. China}
\email{holderhandsome@bjtu.edu.cn}

\author{Yan-Quan Feng}
\address{Yan-Quan Feng, Department of Mathematics, Beijing Jiaotong University, Beijing,
100044, P.R. China}
\email{yqfeng@bjtu.edu.cn}

\author{Dimitri Leemans}
\address{Dimitri Leemans, D\'epartement de Math\'ematique, Universit\'e Libre de Bruxelles, C.P.216, Boulevard du Triomphe, 1050 Bruxelles Belgium}
\email{leemans.dimitri@ulb.be}

\date{}
\maketitle

\begin{abstract}

In [Problems on polytopes, their groups, and realizations, Periodica Math. Hungarica 53 (2006) 231-255] Schulte and Weiss
proposed the following problem: {\em Characterize regular polytopes of orders $2^np$ for $n$ a positive integer and $p$ an odd prime}. In this paper, we first prove that if a $3$-polytope of order $2^np$ has Schl\"afli type $\{k_1, k_2\}$, then $p \mid k_1$ or $p \mid k_2$. This leads to two classes, up to 
duality, for the Schl\"afli type, namely Type (1) where $k_1=2^sp$ and $k_2=2^t$ and Type (2) 
where $k_1=2^sp$ and $k_2=2^tp$. 
We then show that there exists a regular $3$-polytope of order $2^np$ with Type (1) when $s\geq 2$, $t\geq 2$ and $n\geq s+t+1$ coming from a general construction of regular $3$-polytopes of order $2^n\ell_1\ell_2$ with Schl\"afli type $\{2^s\ell_1,2^t\ell_2\}$ where both $\ell_1$ and $\ell_2$ are odd. 
Furthermore, for $p=3$ and $n \geq 7$, we show that there exists a regular
3-polytope of order $3\cdot2^n$ with type $\{6,2^s\}$ if and only if $2\leq s \leq n-2$ and $s \neq n-3$.
For Type (2), we prove that there exists a regular $3$-polytope of order $2^n\cdot 3$ with Schl\"afli type $\{6, 6\}$ when $n \ge 5$ coming from a general construction of regular $3$-polytopes of Schl\"afli type $\{6,6\}$ with orders $192m^3$, $384m^3$ or $768m^3$, for any positive integer $m$.

\bigskip
\f {\bf Keywords:} Regular $3$-polytope, solvable group, automorphism group.
\end{abstract}

\section{Introduction}
An (abstract) polytope $\P$ is a partially-ordered set endowed with a rank function, satisfying certain conditions that arise naturally from a geometric setting. A
\emph{chain} of $\P$ is a totally ordered subset of $\P$, and a maximal chain is called a \emph{flag} of $\P$. The automorphism group of a polytope always acts freely on flags, and in the case where this action is transitive as well, the polytope is said to be \emph{regular} and the order of the automorphism group is called the \emph{order} of the polytope. The study of abstract regular polytopes has a rich history and been described comprehensively in the book by McMullen and Schulte~\cite{ARP}.

It is a natural question to try to classify all pairs $(\P, G)$, where $\P$ is a regular polytope and $G$ is an automorphism group
acting transitively on the flags of $\P$.
The atlas~\cite{atles1} contains information about all regular polytopes with automorphism group of order at most 2000.
Recently, P. Potočnik {\em et al.}~\cite{GroupOrder6000} computed all regular maps with rotational automorphism group of order at most 6000.

 An interesting case is constituted by the pairs $(\P, G)$ with $G$ simple or almost simple.
The atlas~\cite{atles2} contains all regular polytopes whose automorphism group is an almost simple group $G$ such that $S \leq G\leq \Aut(S)$ and $S$ is a simple group of order less than 1 million,
and more striking results have been obtained for the symmetric groups $S_n$ and alternating groups $A_n$: Fernandes~{\em et al.}~\cite{fl,flm,sympolcorr} classified abstract regular polytopes of ranks $n-1$, $n-2$, $n-3$ and $n-4$ for $S_n$, and
Cameron {\em et al.}~\cite{CFLM2017} showed that the highest rank of an abstract regular polytope for $A_n$ is $\lfloor (n-1)/2 \rfloor$ when $n\geq 12$, which was known to be sharp by Fernandes~{\em et al.}~\cite{flm1,flm2}.
More recently, Cameron et al.~\cite{CFL2024} showed that the number of abstract regular polytopes of rank $n-k$ for $S_{n}$ is 
a constant independent of $n$ when $n$ is at least $2k+3$.
For a prime power $q$, results about the highest possible rank of a regular polytope with given automorphism group were obtained for the linear groups $\PSL(2, q)$ in \cite{psl2q}, $\PGL(2,q)$ in \cite{pgl2q}, $\PSL(3, q)$ and $\PGL(3, q)$ in \cite{psl3q}, $\PSL(4, q)$ in \cite{psl4q}, for the Suzuki simple groups $^2B_2(q) = \Sz(q)$ in \cite{suzuki}, for the small Ree groups $^2G_2(q) = R(q)$ in~\cite{lsv2018} and for some symplectic and orthogonal groups in~\cite{brooksbank2021}. Furthermore, Connor {\em et al.}~\cite{soclepsl2q} classified abstract regular polytopes for almost simple groups $G$ with $\PSL(2, q) \leq G \leq \PGaL(2, q)$. For a broader survey on these kind of results, we refer to~\cite{Lee2019}.

Another interesting case is for the pairs $(\P, G)$ with $G$ solvable.
It is well-known that every solvable group, whose order is not a prime, has a non-trivial normal subgroup.
When it comes to solvable groups, except for abelian groups,
the first ones that come to mind are groups of order $2^n$ or $2^np$, which have proved to be more difficult to understand than others.
Schulte and Weiss~\cite{Problem} proposed the following problem.

\begin{prob}
Characterize the groups of orders $2^n$ or $2^np$, with $n$ a positive integer and $p$ an odd prime, which are automorphism groups of regular or chiral polytopes?
\end{prob}

Conder~\cite{SmallestPolytopes} constructed an infinite family of regular polytopes of type $\{4, 4, \cdots, 4\}$ of order $2^n$, where $n$ depends on the rank, and Cunningham and Pellicer~\cite{GD2016} obtained the classification of regular $3$-polytopes of order $2^n$ and Schl\"afli type $\{4, 2^{n-3}\}$. This classification was also obtained by Loyola~\cite{Loyola} by using the classification of $2$-groups with a cyclic subgroup of order $2^{n-3}$.
Gomi, Loyola and De Las Pe\~{n}as~\cite{sc1024} determined the abstract regular polytopes of order $1024$. The authors~\cite{HFL, HFL1} constructed a regular polytope of order $2^n$ for all possible Schl\"afli types given by Conder~\cite{SmallestPolytopes}, and classified the regular $3$-polytopes of order $2^n$ with Schl\"afli types $\{4, 2^{n-4}\}$ and $\{4, 2^{n-5}\}$.

There appears to be few known examples of regular polytopes with solvable automorphism groups of order $2^np$, apart from some of small order and families of regular 3-polytopes arising from
the so-called tight regular $3$-polytopes, namely those with Schl\"afli type $\{k_1, k_2\}$ and order $2k_1k_2$, obtained by Cunningham and Pellicer~\cite{GD2016}.
In this paper, we prove the following theorem (see Section~\ref{Main Results}).

\begin{theorem}\label{maintheorem}
If a regular $3$-polytope of order $2^np$ has type $\{k_1,k_2\}$, then $p \mid k_1$ or $p \mid k_2$. Moreover, up to duality, there are two types of regular $3$-polytopes of order $2^np$, namely
\begin{itemize}
\item Type (1): $k_1=2^sp$ and $k_2=2^t$;
\item Type (2): $k_1=2^sp$ and $k_2=2^tp$.
\end{itemize}
\end{theorem}

In Section~\ref{type1}, we construct regular $3$-polytopes of order $2^n\ell_1\ell_2$ with Schl\"afli type $\{2^s\ell_1,2^t\ell_2\}$ when $s, t \geq 2$, $n\geq s+t+1$ and where both $\ell_1$ and $\ell_2$ are odd.  This shows there are infinitely many pairwise non-isomorphic regular polytopes of order $2^np$ (by taking $\ell_1=p, \ell_2=1$).
Furthermore, for $k_1=p$ or $2p$, the existing results from~\cite{atles1} seem to suggest that polytopes with arbitrary $k_2$ are difficult to obtain for every prime numbers $p$ . We focus on the case $p=3$ and then $k_1=2p=6$, that is, $\mathcal{P}$ has type $\{6, 2^s\}$ and $|\Aut(\mathcal{P})|=3\cdot 2^n$.
Note that if $\mathcal{P}$ has type $\{3, 2^s\}$ then $\Aut(\mathcal{P})$ is a quotient of $[3, \infty]\cong \PGL(2, \mathbb{Z})$, and this case seems much more difficult.
That is the reason why we only consider the case $\{6, 2^s\}$ here. For $n \leq 9$, one may refer to \cite{atles1} and we summarize the results in Table~\ref{table1}.
\begin{table}
\begin{center}
\begin{tabular}{|c|c|c|c|c|c|c|}
\hline 
Type/Order  & $3\cdot 2^4$ & $3\cdot 2^5$ & $3\cdot 2^6$ & $3\cdot 2^7$ & $3\cdot 2^8$ & $3\cdot 2^9$ \\ \hline
$\{6, 2^2\}$     &       Yes      &         Yes        &         Yes        &        Yes         &         Yes       &    Yes                  \\ \hline
$\{6, 2^3\}$    &        $*$      &         Yes         &         Yes       &       Yes          &         Yes        &     Yes                 \\ \hline
$\{6, 2^4\}$    &         $*$     &          $*$        &          Yes      &         {\bf No}         &           Yes       &       Yes               \\ \hline
$\{6, 2^5\}$    &         $*$     &          $*$        &           $*$      &       Yes         &           {\bf No}         &       Yes               \\ \hline
$\{6, 2^6\}$    &         $*$     &          $*$        &           $*$      &         $*$       &            Yes       &        {\bf No}              \\ \hline 
$\{6, 2^7\}$    &         $*$     &          $*$        &            $*$     &         $*$       &            $*$       &          Yes            \\ \hline
\end{tabular}
\end{center}
\caption{Existence of regular polyhedra of certain types and orders.}\label{table1}
\end{table}
These results suggest the following: for $n\geq 7$, there exists a regular
3-polytopes of order $3\cdot2^n$ with type $\{6,2^s\}$ if and only if $2\leq s \leq n-2$ and $s \neq n-3$.
We then prove this is true and we get the following theorem.

\begin{theorem}\label{62n}
For any $n \geq 7$ and $2 \leq s \leq n-2$, there exists a regular
$3$-polytope of order $3\cdot2^n$ with type $\{6,2^s\}$ if and only if $s \neq n-3$.
\end{theorem}

Finally, in Section~\ref{Case2}, we construct a regular $3$-polytope of Schl\"afli type $\{6,6\}$ with order $192m^3$, $384m^3$ or $768m^3$, for any positive integer $m$. In particular, these give a regular $3$-polytope of order $2^np$ for Type~(1) when $s\geq 2$, $t\geq 2$ and $n\geq s+t+1$, and for Type (2) when $p=3$, $s=t=2$ and  $n \ge 5$.

\section{Background results}\label{backgroud}
\subsection{Group theory}

We use  standard notation for group theory, as in~\cite{GroupBook} for example.
In this section we briefly describe some of the specific aspects of group theory that we need.

\smallskip
Let $G$ be any group.   
For a element $g\in G$, we denote by $o(g)$ the order of $g$.
We define the {\em commutator\/} $[x, y]$ of elements $x$ and $y$ of $G$
by $[x, y]=x^{-1}y^{-1}xy$, and then define the {\em derived subgroup} (or {\em commutator subgroup}) of $G$
as the subgroup $G'$ of $G$ generated by all commutators of $G$.
For any non-negative integer $n$, we define the $n^{th}$ derived group of $G$ by setting
$$G^{(0)}=G, \ G^{(1)}=G' \ \hbox{and } \ G^{(n)}=(G^{(n-1)})', \ \hbox{for any} \ n \geq 1.$$

A group $G$ is called {\em solvable\/}  if $G^{(n)}=1$ for some $n$.
This terminology comes from Galois theory, because a polynomial over a field $\mathbb{F}$ is solvable by radicals
if and only if its Galois group over $\mathbb{F}$ is a solvable group.
Every abelian group and every finite $p$-group is solvable, but every non-abelian simple groups is not solvable.
In fact, the smallest non-abelian simple group $A_5$ is also the smallest non-solvable group. 
A subgroup $H$ of a group $G$ is {\em characteristic} if $H^\sigma\leq H$ for every automorphism $\sigma$ of $G$.

The following  proposition is a basic property of commutators and its proof is straightforward.

\begin{prop}\label{commutator}
Let $G$ be a group. For any $x, y, z \in G$, $[xy, z]=[x, z]^y[y, z]$ and $[x, yz]=[x, z][x, y]^z$.
\end{prop}

We need the following two propositions whose proofs are elementary. Hence we state them without proof.

\begin{prop}\label{solvable}
If $N$ is a normal subgroup of the group $G$, such that both $N$ and $G/N$ are solvable, then so is $G$.
\end{prop}

\begin{prop}\label{freeabelian}
Let $G$ be the free abelian group $\mathbb{Z}\oplus \mathbb{Z} \oplus \mathbb{Z}$ of rank $3$, generated by three elements $x, y$ and $z$
subject to the defining relations $[x,y] = [y, z]=[x, z]=1$.
Then for every positive integer $m$, the subgroup $G_m=\lg x^m, y^m, z^m\rg$ is characteristic in $G$, with index $|G:G_m|=m^3$.
\end{prop}

We use the Reidemeister-Schreier theory to produce a defining presentation
for a subgroup $H$ of finite index in a finitely-presented group $G$.
An easily readable reference for this is~\cite[Chapter IV]{Johnson}, but in practice we use its implementation
as the {\tt Rewrite} function in the {\sc Magma} computational software~\cite{BCP97}.
In Section~\ref{type1} and~\ref{Case2}, we heavily rely on some computer computations to deal with small groups and we acknowledge the invaluable help of {\sc Magma}~~\cite{BCP97}.

For a subgroup $H$ of a group $G$, the {\em kernel} $\Core_G(H)$ of $H$ in $G$ is the largest normal subgroup of $G$ contained in $H$. The following result is called {\em Lucchini's Theorem}.

\begin{prop}{\rm \label{core}\cite[Theorem 2.20]{GroupBook}}
Let $A$ be a cyclic proper subgroup of a finite group $G$, and let $K = \Core_{G}(A)$. Then $|A : K| < |G : A|$, and in
particular, if $|A| \geq |G : A|$, then $K > 1$.
\end{prop}

\subsection{String C-groups}

Abstract regular polytopes and string C-groups are the same mathematical objects. The link between these objects may be found for instance in~\cite[Chapter 2]{ARP}.
We take here the viewpoint of string C-groups because it is the easiest and the most efficient one to define abstract regular polytopes.

Let $G$ be a group and let $S=\{\rho_0,\cdots,\rho_{d-1}\}$ be a ordered set of involutions of $G$ that generate $G$.
For $I\subseteq\{0,\cdots,d-1\}$, let $G_I$ denote the group generated by $\{\rho_i:i\in I\}$.
Suppose that
\begin{itemize}
\item[*]  for any $i,j\in \{0, \ldots, d-1\}$ with $|i-j|>1$, $\rho_i$ and $\rho_j$ commute (the \emph{string
property});
\item[*] for any $I,J\subseteq\{0,\cdots,d-1\}$,
$G_I \ \cap \ G_J=G_{I\ \cap \ J}\ \  (\mbox{the \emph{intersection property}})$.
\end{itemize}
Then the pair $(G,S)$ is called a {\em string C-group representation of rank $d$ of $G$} (or {\em string C-group} for short) and the {\em order} of $(G,S)$ is simply the order of $G$. If $(G,S)$  only satisfies the string property, it is called a {\em string group generated by involutions} or \emph{sggi} for short. The intersection property implies that $S$ is a minimal generating set of $G$.  Given a string $C$-group, it is possible to construct an abstract regular polytope as shown in~\cite[Section 2E]{ARP}. Conversely, given an abstract regular polytope $\Gamma$ and one of its flags, we can construct a string C-group representation of $\Aut(\Gamma)$~\cite[Section 2B]{ARP}. 
Therefore, all the results in this paper stated for regular polytopes are also valid for string C-group representations and vice versa.

The \emph{i-faces} of the regular d-polytope
associated with a string C-group $(G,S)$ are the right cosets of the distinguished subgroup $G_i = \lg \r_j \ |\ j \neq i\rg$ for each $i=0,1,\cdots,d-1$, and two faces are incident just when they intersect as cosets. 
In other words, a regular polytope constructed from a string C-group is a coset geometry in the sense of Jacques Tits~\cite{Tits1957}. The {\em Schl\"afli type} of $(G,S)$ is the ordered set $\{p_1,\cdots, p_{d-1}\}$, where $p_i$ is the order of $\r_{i-1}\r_i$. In this paper we assume that each $p_i$ is at least 3 for otherwise the generated group is a direct product of two smaller groups. If that happens, the string C-group (and the corresponding abstract regular polytope) is called {\em degenerate}. The following proposition is related to degenerate string C-groups of rank $3$.

\begin{prop}\label{degenerate}
For $k \geq 2$, let
\begin{itemize}
  \item [$L_1$]$=\lg \r_0, \r_1, \r_2 \ |\ \r_0^2, \r_1^2, \r_2^2, (\r_0\r_1)^{k}, (\r_1\r_2)^{2}, (\r_0\r_2)^2 \rg$,
  \item [$L_2$]$=\lg \r_0, \r_1, \r_2 \ |\ \r_0^2, \r_1^2, \r_2^2, (\r_0\r_1)^{2}, (\r_1\r_2)^{k}, (\r_0\r_2)^2 \rg$.
\end{itemize}
Then $|L_1|=|L_2|=4k$. In particular,
the listed exponents are the true orders of the corresponding elements.
\end{prop}

The proof of Proposition~\ref{degenerate} is straightforward from the fact that  $L_1=\langle\r_0,\r_1\rangle \times \langle \r_2 \rangle \cong D_{2k}\times \mathbb{Z}_2$ and
$L_2=\langle\r_0\rangle\times\langle\r_1,\r_2\rangle\cong \mathbb{Z}_2\times D_{2k}$, where $D_{2k}$ denotes the dihedral group of order $2k$.

Conder~\cite{SmallestPolytopes} obtained a lower bound for the order of a regular polytope.

\begin{prop}{\rm \cite[Theorem 3.2]{SmallestPolytopes}}\label{lowerbound}\label{leastorder}
If $\mathcal{P}$ is a regular $d$-polytope of type $\{k_1, k_2, \ldots, k_{d-1}\}$, then $\mathcal{P}$ has order at least $2k_1k_2 \ldots k_{d-1}$.
\end{prop}

If the lower bound in Proposition~\ref{lowerbound} is attained, $\mathcal{P}$ is called {\em tight}. Cunningham and Pellicer in~\cite{GD2016} 
classified tight regular polyhedra.

A {\em string Coxeter group $[k_1, k_2, \ldots, k_{d-1}]$} is defined as the following group:
\begin{eqnarray*}
\lg \r_0, \r_1, \ldots, \r_{d-1} \ |\ & \r_i^2=1 \ {\rm for} \ 0 \leq i \leq d-1, (\r_i\r_{i+1})^{k_{i+1}}=1 \ {\rm for} \ 0 \leq i \leq d-2,  \\
  & (\r_i\r_j)^2=1 \ {\rm for} \ 0 \leq i < j-1  < d-1\rg.
\end{eqnarray*}

\begin{prop}{\rm \cite[Theorem 5.3]{SmallestPolytopes}}\label{tight}
For every sequence $(k_1, k_2, \ldots, k_{d-1})$ of $d-1$ even integers greater than $2$, there exists a tight regular $d$-polytope $\mathcal{P}$ of order $2k_1k_2 \ldots k_{d-1}$ and type $\{k_1, k_2, \ldots, k_{d-1}\}$. In particular, one can take the string Coxeter group $[k_1, k_2, \ldots, k_{d-1}]$
with standard generators $\r_0, \r_1, \ldots, \r_{d-1}$, and let $\Aut(\mathcal{P})$ be the quotient obtained by adding all relations of the form $[\r_{i}, (\r_{i+1}\r_{i+2})^2]=1$
for $0 \leq i \leq d-3$ and $[(\r_i\r_{i+1})^2, \r_{i+2}]$ for $0 \leq i \leq d-3$.
\end{prop}

The following proposition is called the {\em quotient criterion} for a string C-group.

\begin{prop}{\rm \cite[Section 2E]{ARP}}\label{stringC}
Let $(\G,\{\r_0, \r_1, \r_2\})$ be an sggi, and let $\Lambda = (\lg \s_{0}, \s_{1}, \s_{2}\rg,$ $\{\s_0,\s_1,\s_2\})$  be a string C-group.
If the mapping $\r_j \mapsto \s_j$ for $j=0 ,1, 2$ induces a homomorphism $\pi : \G \rightarrow \Lambda$, which is one-to-one on the subgroup
 $\lg \r_0, \r_1\rg$ or on $\lg \r_1, \r_2\rg$, then $(\G, \{\r_0, \r_1, \r_2\})$ is also a string $C$-group.
\end{prop}

We next state another proposition to determine when some {\em sggi}'s are string C-groups.

\begin{prop}{\rm \cite [Proposition 2.2]{GD2016}}\label{stringCC}
Let $(G,\{\r_0, \r_1, \r_2\})$ be a string C-group. Let $N=\lg (\r_0\r_1)^k \rg$ or $\lg (\r_1\r_2)^k\rg$
for some $k \geq 2$. If $N \unlhd G$, then $(G/N, \{\r_0N, \r_1N, \r_2N\})$ is a string C-group.
\end{prop}

The following proposition gives some famous string C-groups with type $\{4, 4\}$.

\begin{prop}{\rm \cite[Section 8.3]{HW}}\label{type44}
For $b \geq 2$, let
\begin{itemize}
  \item [$M_1$]$=\lg \r_0, \r_1, \r_2 \ |\ \r_0^2, \r_1^2, \r_2^2, (\r_0\r_1)^{4}, (\r_1\r_2)^{4}, (\r_0\r_2)^2, (\r_1\r_0\r_1\r_2)^b\rg$,
  \item [$M_2$]$=\lg \r_0, \r_1, \r_2 \ |\ \r_0^2, \r_1^2, \r_2^2, (\r_0\r_1)^{4}, (\r_1\r_2)^{4}, (\r_0\r_2)^2, (\r_0\r_1\r_2)^{2b} \rg$.
\end{itemize}
Then $M_1\cong (D_b \times D_b)\rtimes  C_2$ and $M_2\cong (D_b \times D_b)\rtimes (C_2\times C_2)$ wtih $|M_1|=8b^2$ and $|M_2|=16b^2$. The listed exponents are the true orders of the corresponding elements.
\end{prop}

\section{Possible types of regular $3$-polytopes of order $2^np$}\label{Main Results}

In this section, we prove Theorem~\ref{maintheorem}. In order to do so, we determine the types of regular $3$-polytopes of order $2^np$. 

\begin{theorem}\label{nonexists}
For an  odd prime $p$ and positive integers $n$, $k_1$ and $k_2$, let $G$ be a group of order $2^np$ generated by $\r_0, \r_1, \r_2$ such that $\r_0^2=\r_1^2= \r_2^2=(\r_0\r_2)^2=(\r_0\r_1)^{k_1}=(\r_1\r_2)^{k_2}=1$. Then $p \di k_1$ or $p\di k_2$.
\end{theorem}

\begin{proof} Note that $|G|=2^np$. We use induction on $n$.

Assume $n=1$. Then $|G|=2p$. Since $\r_0^2=\r_2^2=(\r_0\r_2)^2=1$, we have $\r_0\r_2=\r_2\r_0$. If $|\langle \r_0,\r_1\rangle|\leq 2$ and $|\langle \r_1,\r_2\rangle|\leq 2$, then $\r_0\r_1=\r_1\r_0$ and $\r_1\r_2=\r_2\r_1$. It follows that $G=\langle \r_0,\r_1,\r_2\rangle$ is abelian, and since $\r_0^2=\r_1^2=\r_2^2=1$, $G$ is a $2$-group, contradicting $|G|=2p$. Thus, either $|\langle \r_0,\r_1\rangle|\geq 3$ or $|\langle \r_1,\r_2\rangle|\geq 3$. If $|\langle \r_0,\r_1\rangle|\geq 3$, then $\r_0$ and $\r_1$ are involutions and $\langle \r_0,\r_1\rangle$ is a dihedral group of order $2o(\r_0\r_1)$. In particular, $2o(\r_0\r_1)\di |G|$, and since $|G|=2p$, we have $o(\r_0\r_1)=p$. This implies that $p\di k_1$ because $(\r_0\r_1)^{k_1}=1$. Similarly, if $|\langle \r_1,\r_2\rangle|\geq 3$ then $p\di k_2$. It follows that $p\di k_1$ or $p\di k_2$, as required.

Assume $n\geq 2$. Let $P$ be a Sylow $p$-subgroup of $G$. Then $|P|=p$.

\medskip 
\f {\bf Claim:} $G$ has a non-trivial normal subgroup of order $2$-power.

Suppose on the contrary that $G$ has no non-trivial normal subgroup of order a power of 2. Let $K$ be a minimal normal subgroup of $G$. Since $|G|=2^np$, the Burnside $p^aq^b$-Theorem implies that $G$ is solvable, and so $K$ is elementary abelian. Since $K$ cannot be a non-trivial $2$-group, we have $K=P$ and hence $P\unlhd G$.

Let $C=C_G(P)$ be the centralizer of $P$ in $G$. Then $C\lhd N_{G}(P)=G$. Let $C_2$ be a Sylow $2$-subgroup of $C$. Then $C=C_2\times P$. This implies that $C_2$ is characteristic in $C$, and since $C\lhd G$, we have $C_2\lhd G$. Since $G$ has no non-trivial normal $2$-subgroup, we have $C_2=1$, that is, $C=P$.

By the $N/C$ Theorem, $N_{G}(P)/C_{G}(P)=G/P \lesssim \Aut(P) \cong Z_{p-1}$, and hence $G/P$ is a cyclic group of order $2^n$. Furthermore, $G/P=\langle \r_0P,\r_1P,\r_2P\rangle$, and $G/P$ has a unique cyclic subgroup of order $2$, say $T/P$. Since $(\r_0P)^2=(\r_1P)^2=(\r_2P)^2=P$, we have that $\r_0P\in T/P$, $\r_1P\in T/P$ and $\r_2P\in T/P$, implying $G/P=\langle \r_0P,\r_1P,\r_2P\rangle=T/P$. It follows that $2^n=2$, contradicting $n\geq 2$. This completes the proof of our Claim.

\medskip
We may now use the fact that $G$ must have a normal subgroup $N$ of order $2^t$ with $t\geq 1$. Clearly, $t\leq n$ as $|G|=2^np$. If $t=n$ then $N$ is the unique normal Sylow $2$-subgroup of $G$, and hence $\langle \r_0,\r_1,\r_2 \rangle \leq N$, which implies that $G=\langle \r_0,\r_1,\r_2\rangle=N$, contradicting $|G|=2^np$. Thus, $t<n$, and $G/N$ is a group of order $2^{n-t}p$, generated by $\r_0N$, $\r_1N$ and $\r_2N$ with $(\r_0N)^2=(\r_1N)^2= (\r_2N)^2=(\r_0\r_2N)^2=(\r_0\r_1N)^{k_1}=(\r_1\r_2N)^{k_2}=1$. By inductive hypothesis, $p\di k_1$ or $p\di k_2$, as required. 
\end{proof}

\begin{proof}[Proof of Theorem~\ref{maintheorem}] Let $p$ be a prime and $n$ a positive integer. Let $\mathcal{P}$ be a regular 3-polytope of order $2^np$ with Schl\"afli type $\{k_1, k_2\}$. Then $\Aut(\mathcal{P})$ is a string C-group of rank $3$ generated by three involutions, say $\r_0$, $\r_1$ and $\r_2$, which satisfy all the relations in Theorem~\ref{nonexists}. Thus, $p\di k_1$ or $p\di k_2$. It follows that up to duality, types of $\mathcal{P}$ consist of two classes: Type (1): $k_1=2^sp$ and $k_2=2^t$, and Type (2): $k_1=2^sp$ and $k_2=2^tp$.
\end{proof}

\section{Regular 3-polytopes for Type (1)}\label{type1}

\subsection{Regular 3-polytopes of type $\{2^sp, 2^t\}$ for $s, t \geq 2$}\label{4.1}

In this section we construct infinitely many regular $3$-polytopes of order $2^np$ with Type~(1), namely, Schl\"afli type $\{2^sp, 2^t\}$.

\begin{theorem}\label{maintheorem2}  Let $s, t \geq 2$, $n\geq s+t+1$, and let $\ell_1$ and $\ell_2$ both be odd. Define
$$G=\left\{
\begin{array}{ll}
\lg \r_0, \r_1, \r_2 \ |\ R(\r_0, \r_1, \r_2), [(\r_0\r_1)^2, \r_2]^{2^{\frac{n-s-t-1}{2}}}\rg, & n-s-t\mbox{ odd }\\
\lg \r_0, \r_1, \r_2 \ |\ R(\r_0, \r_1, \r_2), [(\r_0\r_1)^2, (\r_1\r_2)^2]^{2^{\frac{n-s-t-2}{2}}}\rg, & n-s-t\mbox{ even,}
\end{array}
\right.$$
where 
$$R(\r_0, \r_1, \r_2)= \{\r_0^2, \r_1^2, \r_2^2, (\r_0\r_1)^{2^{s}\ell_1}, (\r_1\r_2)^{2^{t}\ell_2}, (\r_0\r_2)^2, [(\r_0\r_1)^4, \r_2],$$ $$[\r_0,(\r_1\r_2)^4]\}.$$

The group $G$ is solvable and the pair $(G,\{\r_0,\r_1,\r_2\})$ is a string C-group of order $2^{n}\ell_1\ell_2$ with type $\{2^s\ell_1, 2^t\ell_2\}$.
\end{theorem}

\begin{proof}
Write $k_1=2^s\ell_1$ and $k_2=2^t\ell_2$. Let
$$
\begin{array}{ll}
&L_1=\lg \r_0, \r_1, \r_2 \ |\ \r_0^2, \r_1^2, \r_2^2, (\r_0\r_1)^{2^s\ell_1}, (\r_1\r_2)^{2}, (\r_0\r_2)^2\rg,\\
&L_2=\lg \r_0, \r_1, \r_2 \ |\ \r_0^2, \r_1^2, \r_2^2, (\r_0\r_1)^2, (\r_1\r_2)^{2^{t}\ell_2}, (\r_0\r_2)^2\rg.\\
\end{array}
$$
By Proposition~\ref{degenerate}, $L_1=\langle\r_0,\r_1\rangle \times \langle \r_2 \rangle \cong D_{2k_1}\times \mathbb{Z}_2$ and
$L_2=\langle\r_0\rangle\times\langle\r_1,\r_2\rangle\cong \mathbb{Z}_2\times D_{2k_2}$.

Since $[\r_0,\r_2]=1$, Proposition~\ref{commutator} gives  $$[(\r_0\r_1)^2,\r_2]=[\r_0\r_1\r_0\r_1,\r_2]=[\r_1\r_0\r_1,\r_2]=[\r_0,\r_1\r_2
\r_1\r_2\r_2]^{\r_1}$$ $$=[\r_0,(\r_1\r_2)^2]^{\r_2\r_1}.$$
The generators $\r_0, \r_1, \r_2$ in $L_1$
satisfy all relations in $G$, and hence $L_1$ is an epimorphic image of $G$ induced by $\r_i\mapsto \r_i$ for $0\leq i\leq 2$. Since $\r_0\r_1$ has order $2^s\ell_1$ in $L_1$, we have that $o(\r_0\r_1)=2^s\ell_1$ in $G$, because $(\r_0\r_1)^{2^{s}\ell_1}=1$ in $G$. Write $A=\lg (\r_0\r_1)^4\rg$. Then $A\leq G$ and   $|A|=o((\r_0\r_1)^4)=(2^s\ell_1)/4=2^{s-2}\ell_1$.

Note that $\lg \r_0, \r_1\rg$ is a dihedral subgroup of $G$. Then $\lg \r_0\r_1\rg \unlhd \lg \r_0, \r_1\rg$. Since $\lg (\r_0\r_1)^4 \rg$ is characteristic in $\lg \r_0\r_1\rg$, we have $\lg (\r_0\r_1)^4 \rg \unlhd \lg \r_0, \r_1\rg$, and since $[(\r_0\r_1)^4,\r_2]=1$ in $G$, we have $\lg (\r_0\r_1)^4 \rg \unlhd G$, that is, $A\unlhd G$.
It follows that$|G|=|G/A| \cdot |A|=|G/A|\cdot 2^{s-2}\ell_1$.
Clearly, $G/A\cong G_1$, where
$$G_1=\left\{
\begin{array}{ll}
\lg \r_0, \r_1, \r_2 \ |\ R_1(\r_0, \r_1, \r_2), [(\r_0\r_1)^2, \r_2]^{2^{\frac{n-s-t-1}{2}}}\rg, & n-s-t\mbox{ odd }\\
\lg \r_0, \r_1, \r_2 \ |\ R_1(\r_0, \r_1, \r_2), [(\r_0\r_1)^2, (\r_1\r_2)^2]^{2^{\frac{n-s-t-2}{2}}}\rg, & n-s-t\mbox{ even,}
\end{array}
\right.$$
with $R_1(\r_0, \r_1, \r_2)= \{\r_0^2, \r_1^2, \r_2^2, (\r_0\r_1)^{4}, (\r_1\r_2)^{2^{t}\ell_2}, (\r_0\r_2)^2, [\r_0,(\r_1\r_2)^4]\}$.

The generators $\r_0, \r_1, \r_2$ in $L_2$ satisfy all relations in $G_1$, and hence $L_2$ is an epimorphic image of $G_2$. Since $\r_1\r_2$ has order $2^t\ell_2$ in $L_2$, we have
that $o(\r_1\r_2)=2^t\ell_2$ in $G_1$ and also in $G$. Let $B=\lg (\r_1\r_2)^4\rg$ be a subgroup of $G_1$. Then $|B|=o((\r_1\r_2)^4)=k_2/4=2^{t-2}\ell_2$.
Since $\langle \r_1,\r_2\rangle$ is a dihedral subgroup of $G_1$ and $[\r_0,(\r_1\r_2)^4]=1$, we have $B\unlhd G_1$ and then $|G_1|=|G_1/B| \cdot |B|=|G_1/B|\cdot 2^{t-2}\ell_2$. In particular, $G_1/B\cong G_2$, where
$$G_2=\left\{
\begin{array}{ll}
\lg \r_0, \r_1, \r_2 \ |\ R_2(\r_0, \r_1, \r_2), [(\r_0\r_1)^2, \r_2]^{2^{\frac{n-s-t-1}{2}}}\rg, & n-s-t\mbox{ odd }\\
\lg \r_0, \r_1, \r_2 \ |\ R_2(\r_0, \r_1, \r_2), [(\r_0\r_1)^2, (\r_1\r_2)^2]^{2^{\frac{n-s-t-2}{2}}}\rg, & n-s-t\mbox{ even,}
\end{array}
\right.$$
with $R_2(\r_0, \r_1, \r_2)= \{\r_0^2, \r_1^2, \r_2^2, (\r_0\r_1)^{4}, (\r_1\r_2)^{4}, (\r_0\r_2)^2\}$.

In the group $G_2$, $(\r_0\r_1)^4=(\r_1\r_2)^4=1$, and since $\r_0\r_2=\r_2\r_0$, we have
$$
\begin{array}{ll}
&[(\r_0\r_1)^2, \r_2]=(\r_1\r_0)^2\r_2(\r_0\r_1)^2\r_2=(\r_1\r_0\r_1\r_2)^2,\\
&[(\r_0\r_1)^2, (\r_1\r_2)^2]=((\r_0\r_1)^2(\r_1\r_2)^2)^2=(\r_0\r_1\r_0\r_2\r_1\r_2)^2=(\r_0\r_1\r_2)^4.\\
\end{array}
$$
It follows that
$$
\begin{array}{ll}
&[(\r_0\r_1)^2, \r_2]^{2^{\frac{n-s-t-1}{2}}}=(\r_1\r_0\r_1\r_2)^{2^{\frac{n-s-t+1}{2}}} \mbox{\hskip 0.5cm with $n-s-t$ odd, and }\\
&[(\r_0\r_1)^2, (\r_1\r_2)^2]^{2^{\frac{n-s-t-2}{2}}}=(\r_0\r_1\r_2){2\cdot 2^{\frac{n-s-t}{2}}} \mbox{\hskip 0.5cm with $n-s-t$ even.}
\end{array}
$$
By Proposition~\ref{type44}, $(G_2, \{\r_0, \r_1, \r_2\})$ is a string C-group of order $2^{n-s-t+4}$ with type $\{4, 4\}$. It follows that $|G_1|=|G_1/B|\cdot 2^{t-2}\ell_2=|G_2|\cdot 2^{t-2}\ell_2=2^{n-s+2}\ell_2$ and
$|G|=|G/A|\cdot 2^{s-2}\ell_1=|G_1|\cdot 2^{s-2}\ell_1=2^n\ell_1\ell_2$.

The generators $\r_0$, $\r_1$, $\r_2$ of $G_2$ satisfy the relations in $G_1$, and hence there is an  epimorphism $\pi: G_1\mapsto G_2$, induced by $\r_0\mapsto \r_0$, $\r_1\mapsto \r_1$ and $\r_2\mapsto \r_2$. Note that $o(\r_0\r_1)=4$ and $o(\r_0)=o(\r_1)=o(\r_2)=o(\r_0\r_2)=2$ in $G_2$. Then $o(\r_0\r_1)=4$ and  $o(\r_0)=o(\r_1)=o(\r_2)=o(\r_0\r_2)=2$ in $G_1$. This implies that $G_1$ is an sggi and the restriction of $\pi$ on the dihedral subgroup $\langle \r_0,\r_1\rangle$ of $G_1$ is one-to-one. By Proposition~\ref{stringC}, $(G_1, \{\r_0, \r_1, \r_2\})$ is a string C-group with type $\{4,k_2\}$.

Also the generators $\r_0$, $\r_1$, $\r_2$ of $G_1$ satisfy the relations in $G$, and hence there is an  epimorphism $G\mapsto G_1$, induced by $\r_0\mapsto \r_0$, $\r_1\mapsto \r_1$ and $\r_2\mapsto \r_2$. We can now use an argument similar to the one above to show that $(G, \{\r_0, \r_1, \r_2\})$ is a string C-group with type $\{k_1, k_2\}$.

By Proposition~\ref{type44}, $G_2 \cong (D_{2^{\frac{n-s-t}{2}}} \times D_{2^{\frac{n-s-t}{2}}})\rtimes (C_2\times C_2)$ or $(D_{2^{\frac{n-s-t+1}{2}}} \times D_{2^{\frac{n-s-t+1}{2}}})\rtimes C_2$, and so $G_2$ is solvable.
Since $G_2\cong G_1/B$ and $B$ is cyclic, $G_1$ is solvable by Proposition~\ref{solvable}, and since $G_1\cong G/A$ and $A$ is cyclic, $G$ is solvable.
\end{proof}

The proof of the above theorem highlights part of the structure of $G$. Indeed, $G$ is a cyclic extension of $G_1$ which itself is a cyclic extension of $G_2$.

Let $\mathcal{P}$ is a regular $3$-polytope of order $2^np$ with Schl\"afli type $\{k_1,k_2\}$, where $k_1=2^s\ell_1$ and $k_2=2^t\ell_2$ with $\ell_1, \ell_2$ odd. By Conder~\cite{SmallestPolytopes}, $|\Aut(\mathcal{P})| \geq 2k_1k_2$, and hence $n\geq s+t+1$. By taking $\ell_1=p$ and $\ell_2=1$ in Theorem~\ref{maintheorem2}, we get the following corollary.

\begin{cor}\label{maintheorem1}
For any odd prime $p$ and positive integers $n, s, t$ with $s\geq 2$, $t\geq 2$ and $n\geq s+t+1$, there exists a regular $3$-polytope of order $2^np$ with Schl\"afli type $\{2^sp, 2^t\}$.
\end{cor}

\subsection{Existence of regular 3-polytopes of type $\{6, 2^s\}$ and order $3\cdot 2^n$}\label{4.2}

In this section, we prove Theorem~\ref{62n}, that is, we prove that there exists a regular
3-polytopes of order $3\cdot2^n$ with type $\{6,2^s\}$ if and only if $2\leq s \leq n-2$ and $s \neq n-3$, where $n \geq 7$.
For $n \leq 9$, one may refer to \cite{atles1}.

 To prove Theorem~\ref{62n}, we first construct infinitely many regular $3$-polytopes of order $96m^3, 192m^3$ and $384m^3$ of type $\{6,4\}$ for $m \geq 1$. Let us start by a general construction. 
For $m\geq 1$, set

\begin{eqnarray*}
R &=& \{\r_0^2, \r_1^2, \r_2^2, (\r_0\r_1)^6, (\r_1\r_2)^4, (\r_0\r_2)^2, [\r_0, (\r_1\r_2)^2]^2\},\\
  R_1 &=& \{(\r_1\r_2^{\r_1\r_0})^{2m}, (\r_2^{\r_1}\r_1^{\r_0})^{2m}, (\r_2\r_1^{\r_0\r_1})^{2m}\}, \\ \vspace{3ex}
  R_2 &=& \{(\r_0^{\r_1}\r_1^{\r_0\r_2})^{2m}, (\r_1^{\r_0\r_2\r_1}\r_0)^{2m}, (\r_0^{\r_1\r_2\r_1}\r_1^{\r_0\r_1})^{2m}\}, \\
  R_3&=& \{(\r_1\r_0\r_2)^{6m}, (\r_2\r_1\r_0)^{6m}, ((\r_2\r_1\r_0)^{\r_1})^{6m}\}.
\end{eqnarray*}

Define
 $$G_m=\lg\, \r_0, \r_1, \r_2 \ | \ R, R_1 \rg,$$
$$H_m=\lg\, \r_0, \r_1, \r_2 \ | \ R, R_2 \rg,$$
$$I_m=\lg\, \r_0, \r_1, \r_2 \ | \ R, R_3 \rg.$$

\begin{theorem}\label{maintheorem3}
$(G_m,\{\r_0,\r_1,\r_2\})$, $(H_m,\{\r_0,\r_1,\r_2\})$ and $(I_m,\{\r_0,\r_1,\r_2\})$ are  solvable string C-groups of type $\{6,4\}$ and respective orders $96m^3, 192m^3$ and $384m^3$.
\end{theorem}

\begin{proof}
We begin by defining ${\mathcal U}$ as the following finitely-presented group.
$$\mathcal{U}:= \lg\, \r_0, \r_1, \r_2 \ | \ \r_0^2, \r_1^2, \r_2^2, (\r_0\r_1)^6, (\r_1\r_2)^4, (\r_0\r_2)^2, [\r_0, (\r_1\r_2)^2]^2 \rg.$$

Using {\sc Magma}~\cite{BCP97}, we can check that this group ${\mathcal U}$ has four normal subgroups of respective index $96, 192$ twice and $384$. Out of the two of index 192, only one can be used to construct an infinite family, namely the one that is a free abelian group of rank three. Hence we have three subgroups at our disposal, namely the subgroups
generated by $\{x_1, y_1, z_1\}, \{x_2, y_2, z_2\}$ and $\{x_3, y_3, z_3\}$, where\\[-2pt]
\begin{center}
\begin{tabular}{lll}
$x_1 = (\r_1\r_2^{\r_1\r_0})^{2}$, \qquad & $x_2=(\r_0^{\r_1}\r_1^{\r_0\r_2})^{2}$ \quad &  and \quad $x_3 = (\r_1\r_0\r_2)^{6}$;\\
$y_1 = (\r_2^{\r_1}\r_1^{\r_0})^{2}$, \qquad & $y_2=(\r_1^{\r_0\r_2\r_1}\r_0)^{2}$ \quad &  and \quad $y_3 =  (\r_2\r_1\r_0)^{6}$;  \\
$z_1 =  (\r_2\r_1^{\r_0\r_1})^{2}$, \qquad & $z_2=(\r_0^{\r_1\r_2\r_1}\r_1^{\r_0\r_1})^{2}$ \quad &  and \quad $z_3 = ((\r_2\r_1\r_0)^{\r_1})^{6}$.
\end{tabular}
\end{center}
The quotients of ${\mathcal U}$ by each of these give the initial members of our three infinite families. Take $N=\lg x_1,y_1,z_1\rg$, $L=\lg x_2,y_2,z_2\rg$ and $M=\lg x_3,y_3,z_3\rg$. Then $N,L,M$ are subgroups of ${\mathcal U}$.

\smallskip
A short computation with {\sc Magma} shows that $N$ is normal in ${\mathcal U}$, with index $96$,
and the {\tt Rewrite} command gives a defining presentation for $N$ with three relators, namely $[x_1, y_1] = [x_1, z_1] = [y_1, z_1]= 1$. It follows that $N\cong \mz\times\mz\times\mz$.
The quotient ${\mathcal U}/N$ is isomorphic to the automorphism group of the regular 3-polytope
of type $\{6, 4\}$ with 96 automorphisms listed at~\cite{atles1}. Furthermore, $({\mathcal U}/N,\{\r_0N,\r_1N,\r_2N\})$ is a string C-group, and $\langle \r_0N,\r_1N\rangle$ is a subgroup of order $12$ in ${\mathcal U}/N$. Similarly, $L$ and $M$ are also free abelian normal groups of rank $3$, $({\mathcal U}/L,\{\r_0L,\r_1L,\r_2L\})$ and $({\mathcal U}/M,\{\r_0M,\r_1M,\r_2M\})$ are string C-groups of order $192$ and $384$ respectively, and $|\langle \r_0L,\r_1L\rangle|=|\langle \r_0M,\r_1M\rangle|=12$.

Take $N_m=\lg x_1^m, y_1^m,z_1^m\rg$, $L_m=\lg x_2^m, y_2^m,z_2^m\rg$ and $M_m=\lg x_3^m, y_3^m,z_3^m\rg$. Then $N_1=N$, $L_1=L$ and $M_1=M$. By Proposition~\ref{freeabelian}, $N_m$, $L_m$ and $M_m$ are characteristic in $N$ and hence normal in ${\mathcal U}$, with index $|\,{\mathcal U}:N_m| = |\,{\mathcal U}:N||N:N_m| = 96m^3$, $|\,{\mathcal U}:L_m| = |\,{\mathcal U}:L||L:L_m| =192m^3$ and $|\,{\mathcal U}:M_m| = |\,{\mathcal U}:M||M:M_m| = 384m^3$.

The subgroup  $N/N_m$ of ${\mathcal U}/N_m$ is abelian and normal. Also $({\mathcal U}/N_m)/(N/N_m) \cong  {\mathcal U}/N$ is a $\{2, 3\}$-group of order $96$. By Proposition~\ref{solvable}, ${\mathcal U}/N_m$ is solvable. Similarly, ${\mathcal U}/L_m$ and ${\mathcal U}/M_m$ are solvable.

Since $N_m\leq N$, the map $\r_0N_m\mapsto \r_0N$, $\r_1N_m\mapsto \r_1N$ and $\r_2N_m\mapsto \r_2N$ induces an epimorphism from ${\mathcal U}/N_m$ to ${\mathcal U}/N$, say $\a$. Noting that  $(\r_0N_m)^2=(\r_1N_m)^2=(\r_0\r_1N_m)^6=1$, we have $|\langle \r_0N_m,\r_1N_m\rangle|\leq 12$ in ${\mathcal U}/N_m$, and hence $|\langle \r_0N_m,\r_1N_m\rangle|=12$ because $|\langle \r_0N,\r_1N\rangle|=12$. This implies that the restriction of $\a$ on $\lg \r_0N_m, \r_1N_m \rg$ is a bijection from $\lg \r_0N_m, \r_1N_m \rg$ to $\lg \r_0N, \r_1N \rg$, and by  Proposition~\ref{stringC}, $({\mathcal U}/N_m,\{\r_0N_m,\r_1N_m,\r_2N_m\})$ is a string C-group of type $\{6, 4\}$ and order $96m^3$. 

Similarly, $({\mathcal U}/L_m,\{\r_0L_m,\r_1L_m,\r_2L_m\})$ and $({\mathcal U}/M_m,\{\r_0M_m,\r_1M_m,\r_2M_m\})$ are  solvable string C-groups of type $\{6,4\}$ and respective orders $192m^3$ and $384m^3$. Clearly, $G_m\cong {\mathcal U}/N_m$, $H_m\cong {\mathcal U}/L_m$ and $I_m\cong {\mathcal U}/M_m$. This completes the proof.  \end{proof}

By taking $m$ a power of $2$ in the previous theorem,
we have the following corollary.

\begin{cor}\label{4.5}
For any $n \geq 5$, let $$R= \{\r_0^2, \r_1^2, \r_2^2, (\r_0\r_1)^{6}, (\r_1\r_2)^{4}, (\r_0\r_2)^{2}, [\r_0,(\r_1\r_2)^2]^2\},$$

and for $n \equiv i \ \mod \ 3$ let

$$G := \lg \r_0, \r_1, \r_2 \ |\ R, R_i\rg$$

where 

$$R_0 := \{(\r_0^{\r_1}\r_1^{\r_0\r_2})^{2\cdot 2^{\frac{n-6}{3}}}, (\r_1^{\r_0\r_2\r_1}\r_0)^{2\cdot 2^{\frac{n-6}{3}}}, (\r_0^{\r_1\r_2\r_1}\r_1^{\r_0\r_1})^{2\cdot 2^{\frac{n-6}{3}}}\},$$

$$R_1:=\{ (\r_1\r_0\r_2)^{6\cdot 2^{\frac{n-7}{3}}}, (\r_2\r_1\r_0)^{6\cdot 2^{\frac{n-7}{3}}}, ((\r_2\r_1\r_0)^{\r_1})^{6\cdot 2^{\frac{n-7}{3}}}\},$$

$$R_2:=\{(\r_1\r_2^{\r_1\r_0})^{2\cdot 2^{\frac{n-5}{3}}}, (\r_2^{\r_1}\r_1^{\r_0})^{2\cdot 2^{\frac{n-5}{3}}}, (\r_2\r_1^{\r_0\r_1})^{2\cdot 2^{\frac{n-5}{3}}}\}.$$

Then $G$ is the automorphism group of a regular $3$-polytope with type $\{6, 4\}$ and order $3\cdot2^n$. 
\end{cor}

\begin{lem}\label{anytype}
For any $n \geq 7$, let $x=\r_1\r_2$,
\begin{center}
\begin{tabular}{lll}
$x_1 = (\r_1\r_2^{\r_1\r_0})^{2}$, \qquad & $x_2=(\r_0^{\r_1}\r_1^{\r_0\r_2})^{2}$ \quad &  and \quad $x_3 = (\r_1\r_0\r_2)^{6}$,  \\[+2pt]
$y_1 = (\r_2^{\r_1}\r_1^{\r_0})^{2}$, \qquad & $y_2=(\r_1^{\r_0\r_2\r_1}\r_0)^{2}$ \quad &  and \quad $y_3 =  (\r_2\r_1\r_0)^{6}$,  \\[+2pt]
$z_1 =  (\r_2\r_1^{\r_0\r_1})^{2}$, \qquad & $z_2=(\r_0^{\r_1\r_2\r_1}\r_1^{\r_0\r_1})^{2}$ \quad &  and \quad $z_3 = ((\r_2\r_1\r_0)^{\r_1})^{6}$,  \\[+2pt]
\end{tabular}
\end{center}
$$R'=\{\r_0^2, \r_1^2, \r_2^2, (\r_0\r_1)^{6}, (\r_1\r_2)^{2^s}, (\r_0\r_2)^{2}, [\r_0,(\r_1\r_2)^2]^2, [\r_0, (\r_1\r_2)^4]\}.$$
Let $2\leq s\leq n-4$ and $s\neq n-5$.
For $n -s \equiv i \ \mod \ 3$ let

$$G := \lg \r_0, \r_1, \r_2 \ |\ R', R'_i\rg$$

where 

$$ R'_0 = \{x_1^{2^{\frac{n-s-3}{3}}}=x^{2\cdot 2^{\frac{n-s-3}{3}}}, y_1^{2^{\frac{n-s-3}{3}}}=x^{2\cdot 2^{\frac{n-s-3}{3}}}, z_1^{2^{\frac{n-s-3}{3}}}=(x^{-1})^{2\cdot 2^{\frac{n-s-3}{3}}}\},$$

$$ R'_1 = \{x_2^{ 2^{\frac{n-s-4}{3}}}, y_2^{ 2^{\frac{n-s-4}{3}}}, z_2^{2^{\frac{n-s-4}{3}}}\},$$

$$ R'_2 = \{x_3^{ 2^{\frac{n-s-5}{3}}}=x^{6\cdot 2^{\frac{n-s-5}{3}}}, y_3^{2^{\frac{n-s-5}{3}}}=(x^{-1})^{6\cdot 2^{\frac{n-s-5}{3}}}, z_3^{2^{\frac{n-s-5}{3}}}=x^{6\cdot 2^{\frac{n-s-5}{3}}}\}.$$
Then $G$ is the automorphism group of a regular $3$-polytope with type $\{6, 2^s\}$ and has order $3\cdot2^n$ provided $n-s \neq 5$. 
\end{lem}

\begin{proof}
Let $L=\lg \r_0, \r_1, \r_2 \ |\ \r_0^2, \r_1^2, \r_2^2, (\r_0\r_1)^{2}, (\r_1\r_2)^{2^s}, (\r_0\r_2)^2 \rg \cong \mathbb{Z}_2 \times D_{2^{s+1}}$. It is easy to see that the 
generators $\r_0, \r_1, \r_2$ in $L$ satisfy all the relations given in the definition of $G$. This implies that $L$ is a homomorphic image of $G$.
By Proposition~\ref{degenerate}, $\r_1\r_2$ has order $2^s$ in $L$ and hence has order $2^s$ in $G$. Note that $[\r_0, (\r_1\r_2)^4]=1$
in $G$ implies that $\lg (\r_1\r_2)^4\rg \unlhd G$. 
Observe that if $n-s=5$ the quotient is not smooth because of the relations in $R'_2$, hence the condition that $n-s\neq 5$ in the statement of the lemma.
It follows that $|G|=o((\r_1\r_2)^4)\cdot |G/\lg (\r_1\r_2)^4\rg|=2^{s-2}\cdot |G/\lg (\r_1\r_2)^4\rg|$.
Write $\olg=G/\lg (\r_1\r_2)^4\rg$, where $\overline{g}=g\lg (\r_1\r_2)^4\rg$ for any $g \in G$.

Now, $\olg$ has the following presentation
\begin{small}
$$\olg=\left\{
\begin{array}{ll}
\lg \overline{\r_0}, \overline{\r_1}, \overline{\r_2} \ |\ \overline{R'}, \overline{y}^{2^{\frac{n-s-4}{3}}}\rg, & n-s \equiv 1 \ \mod 3\\
\lg \overline{\r_0}, \overline{\r_1}, \overline{\r_2} \ |\ \overline{R'}, \overline{y}^{2^{\frac{n-s-2}{3}}}, (\overline{\r_0\r_1\r_2})^{3\cdot 2^{\frac{n-s-2}{3}}}\rg, & n-s \equiv 2 \ \mod 3\\
\lg \overline{\r_0}, \overline{\r_1}, \overline{\r_2} \ |\ \overline{R'}, \overline{y}^{2^{\frac{n-s-3}{3}}}, [(\overline{(\r_0\r_1)^3(\r_1\r_2)^2})^{2^{\frac{n-s-3}{3}}}, \overline{\r_2}]\rg, & n-s \equiv 0  \ \mod 3,
\end{array}
\right.$$
\end{small}
where $n-4 \geq s \geq 2, n-s \neq 5$ and 
$$\overline{R}=\overline{R}(\overline{\r_0}, \overline{\r_1}, \overline{\r_2})= \{\overline{\r_0^2}, \overline{\r_1^2}, \overline{\r_2^2}, \overline{(\r_0\r_1)^{6}}, \overline{(\r_1\r_2)^{4}}, \overline{(\r_0\r_2)^{2}}, [\overline{\r_0},(\overline{\r_1\r_2})^2]^2\}.$$
By Corollary~\ref{4.5}, we have $|\olg|=3\cdot 2^{n-s+2}$, and hence $|G|=3\cdot 2^n$.

Since $(\olg, \{\overline{\r_0}, \overline{\r_1}, \overline{\r_2}\})$ is a string C-group and $\r_0\r_1$ has order $6$ in $G$ and $\olg$,
by Proposition~\ref{stringC}, we have that $\{G, (\r_0, \r_1, \r_2)\}$ is a string C-group.
\end{proof}

\begin{lem}\label{type$n-s=5$}
For any $n \geq 7$, let 
$$G=\lg \r_0, \r_1, \r_2 \ |\ \r_0^2, \r_1^2, \r_2^2, (\r_0\r_1)^{6}, (\r_1\r_2)^{2^{n-5}}, (\r_0\r_2)^{2}, [\r_0, (\r_1\r_2)^4],$$ $$(\r_1\r_0\r_2\r_1\r_0(\r_2\r_1)^2\r_0)^2\rg.$$
Then $G$ is the automorphism group of a regular $3$-polytope with type $\{6, 2^{n-5}\}$ and has order $3\cdot2^n$. 
\end{lem}

\begin{proof} Let $L=\lg \r_0, \r_1, \r_2 \ |\ \r_0^2, \r_1^2, \r_2^2, (\r_0\r_1)^{2}, (\r_1\r_2)^{2^{n-5}}, (\r_0\r_2)^2 \rg$. It is easy to see that the 
generators $\r_0, \r_1, \r_2$ in $L$ satisfy all the relations given in the definition of $G$. This implies that $L$ is a homomorphic image of $G$.
By Proposition~\ref{degenerate}, $\r_1\r_2$ has order $2^{n-5}$ in $L$ and hence has order $2^{n-5}$ in $G$. 

Since $[\r_0, (\r_1\r_2)^4]=1$ and $((\r_1\r_2)^4)^{\r_1}=((\r_1\r_2)^4)^{\r_2}=(\r_1\r_2)^{-4}$, we have $\lg (\r_1\r_2)^4\rg \unlhd G$.
Consider the quotient group $\olg=G/\lg (\r_1\r_2)^4\rg$, where $\overline{g}=g\lg (\r_1\r_2)^4\rg$ for any $g \in G$. Then $|G|=2^{n-7}\cdot |\olg|$.
Using {\sc Magma}, we get that $|\olg|=3\cdot2^7$ and $(\olg,\{\overline{\r_0}, \overline{\r_1}, \overline{\r_2}\})$ is a string C-group with type $\{6, 4\}$. Then $|G|=3\cdot2^n$
By Proposition~\ref{stringC}, we have $(G, \{\r_0, \r_1, \r_2\})$ is a string C-group with type $\{6, 2^{n-5}\}$.
\end{proof}
Observe that when $n=7$ in the above lemma, we get a 3-polytope of type $\{6,4\}$ and order $3\cdot 2^7 = 384$. It is not isomorphic to the regular 3-polytope of type $\{6,4\}$ and order 384 given by Corollary~\ref{4.5}. These two are the only two regular 3-polytopes of type $\{6,4\}$ and order 384, as it can be checked in~\cite{atles1}.

\begin{lem}\label{type$n-s=3$}
There is no regular $3$-polytopes of order $3\cdot 2^n$ with type $\{6, 2^{n-3}\}$ for $n\geq 7$.
\end{lem}

\begin{proof}
Using {\sc Magma}~\cite{BCP97}, we can check that that there is no string C-group of order $3\cdot 2^n$ with type $\{6, 2^{n-3}\}$ for $n=7$.

Suppose there exists a string C-group $(G, \{\r_0, \r_1, \r_2\})$ of order $3\cdot 2^{n}$ with type $\{6, 2^{n-3}\}$ for $n \geq 8$.
Since $2^{n-3}\cdot 2^{n-3}=2^{n-6}\cdot 2^{n}\geq 3\cdot 2^n$ for $n \geq 8$, by Proposition~\ref{core}, we have $|$Core$_{G}(\r_1\r_2)|>1$.
It follows that $\lg (\r_1\r_2)^{2^{n-4}}\rg \unlhd G$. Consider the quotient group $\olg=G/\lg (\r_1\r_2)^{2^{n-4}}\rg$ where $\overline{g}=g\lg(\r_1\r_2)^{2^{n-4}}\rg$
for every $g \in G$. By Proposition~\ref{stringCC}, $(\olg, \{\overline{\r_0}, \overline{\r_1}, \overline{\r_2}\})$ is a string C-group of order $3\cdot 2^{n-1}$ with type $\{6, 2^{n-4}\}$.
This, combined with the fact that there is no such string C-group for $n=7$ gives a contradiction. Hence there is no string C-group of the above type and order for $n\geq 8$ as well.
\end{proof}

\begin{proof}[Proof of Theorem~\ref{62n}]
Proposition~\ref{tight} deals with the case $s=n-2$ and shows that there exists a regular 3-polytope of order $3\cdot 2^n$ with type $\{6,2^{n-2}\}$.
Lemma~\ref{type$n-s=3$} shows that there is no regular 3-polytope in the case where $s=n-3$ and $n\geq 7$.
Corollary~\ref{4.5} shows that there exist regular 3-polytopes of order $3\cdot 2^n$ and type $\{6,4\}$ for all $n\geq 5$.
Lemmas~\ref{anytype} and~\ref{type$n-s=5$} show that 
there exist regular 3-polytopes of order $3\cdot 2^n$ and type $\{6,2^s\}$ for all $n\geq 7$ and $s\leq n-4$. The few remaining small cases can be handled with {\sc Magma}.
\end{proof}
\section{Regular 3-polytopes for Type (2)\label{Case2}}

In this section, we construct infinitely many regular $3$-polytopes of order $3\cdot 2^n$ of type $\{6,6\}$. Let us start by a general construction. For $m\geq 1$, set

\begin{eqnarray*}
R&=&\r_0^2, \r_1^2, \r_2^2, (\r_0\r_1)^6, (\r_1\r_2)^6, (\r_0\r_2)^2, (\r_2\r_1\r_0\r_1)^3,\\
  R_1 &=& \{(\r_0\r_2\r_1)^{4m})^{\r_1},(\r_2\r_1\r_0)^{4m},(\r_2\r_1\r_0)^{4m})^{\r_1}\}, \\ \vspace{3ex}
  R_2 &=& \{(((\r_1\r_2)^3(\r_0\r_1)^3)^{2m})^{\r_0\r_1},(((\r_1\r_2)^3(\r_0\r_1)^3)^{2m})^{\r_2}, (((\r_1\r_2)^3(\r_0\r_1)^3)^{2m})^{\r_2\r_1}\}, \\
  R_3&=& \{(((\r_0\r_1)^2(\r_2\r_1)^2)^{3m})^{\r_0},((\r_1\r_2)^2(\r_1\r_0)^2)^{3m}, (((\r_1\r_0)^2(\r_1\r_2)^2)^{3m})^{\r_2}\}.
\end{eqnarray*}

Define
 $$G_m=\lg\, \r_0, \r_1, \r_2 \ | \ R, R_1 \rg,$$
$$H_m=\lg\, \r_0, \r_1, \r_2 \ | \ R, R_2 \rg,$$
$$I_m=\lg\, \r_0, \r_1, \r_2 \ | \ R, R_3 \rg,$$

\begin{theorem}\label{maintheorem4}
$(G_m,\{\r_0,\r_1,\r_2\})$, $(H_m,\{\r_0,\r_1,\r_2\})$ and $(I_m,\{\r_0,\r_1,\r_2\})$ are  solvable string C-groups of type $\{6,6\}$ of order $192m^3, 384m^3$ and $768m^3$, respectively.
\end{theorem}

\demo
We begin by defining ${\mathcal U}$ as the following finitely-presented group.
$$\mathcal{U} := \lg\, \r_0, \r_1, \r_2 \ | \ \r_0^2, \r_1^2, \r_2^2, (\r_0\r_1)^6, (\r_1\r_2)^6, (\r_0\r_2)^2, (\r_2\r_1\r_0\r_1)^3 \,\rg$$

Using {\sc Magma}~\cite{BCP97}, we can see that this group ${\mathcal U}$ has three normal subgroups of index $192, 384$ and $768$, namely the subgroups
generated by $\{x_1, y_1, z_1\}$, $\{x_2, y_2, z_2\}$ and $\{x_3, y_3, z_3\}$, where\\[-2pt]
\begin{center}
\begin{tabular}{lll}
$x_1 = (\r_0\r_2\r_1)^4)^{\r_1}$, \qquad & $x_2= (((\r_1\r_2)^3(\r_0\r_1)^3)^2)^{\r_0\r_1}$, \quad & $x_3 = (((\r_0\r_1)^2(\r_2\r_1)^2)^3)^{\r_0}$,  \\[+2pt]
$y_1 = (\r_2\r_1\r_0)^4$, \qquad & $y_2= (((\r_1\r_2)^3(\r_0\r_1)^3)^2)^{\r_2}$, \quad & $y_3 = ((\r_1\r_2)^2(\r_1\r_0)^2)^3$,  \\[+2pt]
$z_1 = (\r_2\r_1\r_0)^4)^{\r_1}$, \qquad & $z_2= (((\r_1\r_2)^3(\r_0\r_1)^3)^2)^{\r_2\r_1},$ \quad & $z_3 = (((\r_1\r_0)^2(\r_1\r_2)^2)^3)^{\r_2}$,  \\[+2pt]
\end{tabular}
\end{center}
respectively. The quotients of ${\mathcal U}$ by each of these give the initial members of our three infinite families. Write $N=\lg x_1,y_1,z_1\rg$, $L=\lg x_2,y_2,z_2\rg$ and $M=\lg x_3,y_3,z_3\rg$. Then $N,L,M$ are subgroups of ${\mathcal U}$.

\smallskip
A short computation with {\sc Magma} shows that $N$ is normal in ${\mathcal U}$, with index $192$,
and then the {\tt Rewrite} command gives a defining presentation for $N$ with three relations, namely $[x_1, y_1] = [x_1, z_1] = [y_1, z_1]= 1$. It follows that $N\cong \mz\times\mz\times\mz$.
The quotient ${\mathcal U}/N$ is isomorphic to the automorphism group of the regular 3-polytope
of type $\{6, 6\}$ with 192 automorphisms listed at~\cite{atles1}, and so $({\mathcal U}/N,\{\r_0N,\r_1N,\r_2N\})$ is a string C-group of order $192$ with type $\{6, 6\}$. In particular, $\langle \r_0N,\r_1N\rangle$ is a subgroup of ${\mathcal U}/N$ of order $12$. Similarly, $L$ and $M$ are also free abelian normal groups of rank $3$, $({\mathcal U}/L,\{\r_0L,\r_1L,\r_2L\})$ and $({\mathcal U}/M,\{\r_0M,\r_1M,\r_2M\})$ are string C-groups of order $384$ and $768$ respectively, and $|\langle \r_0L,\r_1L\rangle|=|\langle \r_0M,\r_1M\rangle|=12$.

Take $N_m=\lg x_1^m, y_1^m,z_1^m\rg$, $L_m=\lg x_2^m, y_2^m,z_2^m\rg$ and $M_m=\lg x_3^m, y_3^m,z_3^m\rg$. Then $N_1=N$, $L_1=L$ and $M_1=M$. By Proposition~\ref{freeabelian}, $N_m$, $L_m$ and $M_m$ are characteristic in $N$ and hence normal in ${\mathcal U}$, with index $|\,{\mathcal U}:N_m| = |\,{\mathcal U}:N||N:N_m| = 192m^3$, $|\,{\mathcal U}:L_m| = |\,{\mathcal U}:L||L:L_m| =384m^3$ and $|\,{\mathcal U}:M_m| = |\,{\mathcal U}:M||M:M_m| = 768m^3$.

The subgroup  $N/N_m$ of ${\mathcal U}/N_m$ is abelian and normal, and $({\mathcal U}/N_m)/(N/N_m) \cong  {\mathcal U}/N$ is a $\{2, 3\}$-group of order $192$. By Proposition~\ref{solvable}, ${\mathcal U}/N_m$ is solvable. Similarly, ${\mathcal U}/L_m$ and ${\mathcal U}/M_m$ are solvable.

Since $N_m\leq N$, the map $\r_0N_m\mapsto \r_0N$, $\r_1N_m\mapsto \r_1N$ and $\r_2N_m\mapsto \r_2N$ induces an epimorphism from ${\mathcal U}/N_m$ to ${\mathcal U}/N$, say $\a$. Noting that  $(\r_0N_m)^2=(\r_1N_m)^2=(\r_0\r_1N_m)^6=1$, we have $|\langle \r_0N_m,\r_1N_m\rangle|\leq 12$ in ${\mathcal U}/N_m$, and hence $|\langle \r_0N_m,\r_1N_m\rangle|=12$ because $|\langle \r_0N,\r_1N\rangle|=12$. This implies that the restriction of $\a$ on $\lg \r_0N_m, \r_1N_m \rg$ is a bijection from $\lg \r_0N_m, \r_1N_m \rg$ to $\lg \r_0N, \r_1N \rg$, and by  Proposition~\ref{stringC}, $({\mathcal U}/N_m,\{\r_0N_m,\r_1N_m,\r_2N_m\})$ is a string C-group of type $\{6, 6\}$ of order $192m^3$. 

Similarly, $({\mathcal U}/L_m,\{\r_0L_m,\r_1L_m,\r_2L_m\})$ and $({\mathcal U}/M_m,\{\r_0M_m,\r_1M_m,\r_2M_m\})$ are  solvable string C-groups of type $\{6,6\}$ and respective orders $384m^3$ and $768m^3$. Clearly, $G_m\cong {\mathcal U}/N_m$, $H_m\cong {\mathcal U}/L_m$ and $I_m\cong {\mathcal U}/M_m$. This completes the proof.  \hfill\qed

By Conder~\cite{atles1}, there exists a regular $3$-polytope of order $3\cdot 2^5$ with Schl\"afli type $\{6, 6\}$. By taking $m$ as a $2$-power in Theorem~\ref{maintheorem}, we have the following corollary.

\begin{cor}
There exists a regular $3$-polytope of order $3\cdot 2^n$ with Schl\"afli type $\{6, 6\}$, for every integer $n \ge 5$.
\end{cor}

We do not know whether there exists a regular $3$-polytope of order $2^np$ with Type~(2) for every prime $p\geq 5$.

\medskip\medskip

\section{Acknowledgements} This work was supported by the National Natural Science Foundation of China (12201371,
12331013,12311530692,12271024,12161141005), the 111 Project of China (B16002), the Fundamental Research Program of Shanxi Province 20210302124078, an Action de Recherche Concertée grant of the Communauté Française Wallonie Bruxelles and a PINT-BILAT-M grant from the Fonds National de la Recherche Scientifique de Belgique (FRS-FNRS).

\end{document}